\pdfoutput=1
\documentclass[english, compsoc,12pt,final,onecolumn,  semicolon]{IEEEtran}
\usepackage[T1]{fontenc}
\usepackage[latin9]{inputenc}
\usepackage{babel}
\RequirePackage{xcolor}
\usepackage{bm}
\usepackage{amsmath}
\usepackage{amsthm}
\usepackage{amssymb}
\usepackage[numbers,square,sort&compress]{natbib}
 
\RequirePackage{hyperref}
\hypersetup{	
				unicode=true,
 				bookmarks=true,
 				bookmarksnumbered=false,
 				bookmarksopen=false,
 				breaklinks=true,
 				pdfborder={0 0 0},
 				pdfborderstyle={},
 				colorlinks=true,
 				linkcolor=blue!75!green!80!black,
 				citecolor=red!66!black
			}
 			
\makeatletter
\theoremstyle{plain}
\newtheorem{thm}{\protect\theoremname}
\theoremstyle{definition}
\newtheorem{defn}[thm]{\protect\definitionname}
\theoremstyle{plain}
\newtheorem{cor}[thm]{\protect\corollaryname}
\theoremstyle{plain}
\newtheorem{lem}[thm]{\protect\lemmaname}
\theoremstyle{plain}
\newtheorem{prop}[thm]{\protect\propositionname}

\usepackage{microtype, lmodern}
\usepackage[bb=boondox]{mathalfa}
\RequirePackage{graphicx}

\newif\ifonecol
\if@twocolumn
	\onecolfalse
\else
	\onecoltrue
\fi

\makeatother

\providecommand{\corollaryname}{Corollary}
\providecommand{\definitionname}{Definition}
\providecommand{\lemmaname}{Lemma}
\providecommand{\propositionname}{Proposition}
\providecommand{\theoremname}{Theorem}

\global\long\def\mb#1{\hm{#1}}
\global\long\def\mbb#1{\mathbb{#1}}

\global\long\def\mr#1{\mathrm{#1}}
\global\long\def\msf#1{\mathsf{#1}}
\global\long\def\E{\mbb E}
\global\long\def\P{\mbb P}

\global\long\def\T{\scalebox{0.55}{\ensuremath{\msf T}}}

\global\long\def\tr{\mr{tr}}

\global\long\def\defeq{\stackrel{\textup{\tiny def}}{=}}

\begin{document}

\title{Algebraic Connectivity Under Site Percolation\\
in Finite Weighted Graphs}

\author{Sohail~Bahmani, Justin~Romberg,~\IEEEmembership{Senior~Member,~IEEE}, and~Prasad~Tetali\IEEEcompsocitemizethanks{\IEEEcompsocthanksitem S.~Bahmani and J.~Romberg are with the School of Electrical and Computer Engineering,
Georgia Institute of Technology, Atlanta, GA 30332 USA.\protect\\ E-mail: \texttt{\{sohail.bahmani,jrom\}@ece.gatech.edu} \IEEEcompsocthanksitem P.~Tetali is with the School of Mathematics and the College of Computing,
Georgia Institute of Technology, Atlanta, GA 30332 USA.\protect\\ E-mail: \texttt{tetali@math.gatech.edu}}}

\IEEEtitleabstractindextext{%
\begin{abstract}
We study the behavior of \emph{algebraic connectivity} in a weighted
graph that is subject to \emph{site percolation}, random deletion
of the vertices. Using a refined concentration inequality for random
matrices we show in our main theorem that the (augmented) Laplacian
of the percolated graph concentrates around its expectation. This
concentration bound then provides a lower bound on the algebraic connectivity
of the percolated graph. As a special case for $(n,d,\lambda)$-graphs
(i.e., $d$-regular graphs on $n$ vertices with non-trivial eigenvalues
less than $\lambda$ in magnitude) our result shows that, with high
probability, the graph remains connected under a homogeneous site percolation
with survival probability $p\ge1-C_{1}n^{-C_{2}/d}$ with $C_{1}$
and $C_{2}$ depending only on $\lambda/d$.
\end{abstract}

}
\maketitle 
\section{Introduction}

\IEEEPARstart{C}{onsider} a connected weighted graph $G=\left(V=\left[n\right],E\right)$
with (non-negative) edge weights $\left\{ w_{i,j}\right\} _{1\le i,j\le n}$
and no self-loop (i.e., $w_{i,i}=0$ for all $1\le i\le n$) and suppose
that each vertex $i$ of $G$ is deleted independently with probability
$1-p_{i}$. These types of random graph models can describe certain
phenomena in random media and are studied under \emph{percolation
theory} \citep{grimmett_percolation_1999} in mathematics and statistical
physics. The process of vertex deletion, as described above, is usually
referred to as \emph{site percolation} whereas \emph{bond percolation}
refers to the process of random deletion (or addition) of the edges
of a graph. In this paper we establish a lower bound on \emph{algebraic
connectivity} of the surviving subgraph in the described site percolation
model. The algebraic connectivity of a graph $G$ is $a=\lambda_{2}(\mb L)$,
the second smallest eigenvalue of the graph \emph{Laplacian} 
\begin{align*}
\mb L & \defeq\sum_{1\le i<j\le n}w_{i,j}(\mb e_{i}-\mb e_{j})(\mb e_{i}-\mb e_{j})^{\T},
\end{align*}
 where $\mb e_{i}$'s are canonical basis vectors. Algebraic connectivity
and its analog for normalized Laplacians are important because they
provide a bound on isoperimetric constants of graphs through Cheeger's
inequality \citep[see e.g.,][]{mohar_isoperimetric_1989} and they
are critical in approximation of the mixing rate of continuous-time
Markov chains \citep[Ch. 20]{sun_fastest_2006,levin_markov_2009}.

Properties such as connectivity, spectral gap, and emergence of a
\emph{giant component} (i.e., a connected component with $\Omega(n)$
vertices) have received more attention and are better understood for
bond percolation models compared to site percolation models. Perhaps,
the main reason is that in bond percolation edges are removed independently
whereas in site percolation edge deletions are dependent since they
share a common vertex which lead to more intricate behavior.

In this paper we focus on algebraic connectivity of the surviving
subgraph of a weighted graph under (inhomogeneous) site percolation. Using a delicate matrix concentration inequality (Proposition \ref{prp:Bernoulli-series-tail}), in our main result  (Theorem \ref{thm:aug-Laplacian-tail}) we show that the ``augmented'' Laplacian concentrates around its expectation. This result allows us to find a non-trivial lower bound on the algebraic connectivity in a straightforward way. For concreteness, we also apply the general result  of the Theorem \ref{thm:aug-Laplacian-tail} to obtain a threshold for connectivity in the special case of $(n,d,\lambda)$-graphs under uniform site percolation. In particular, Corollary \ref{cor:ndl-graphs} below shows that if the vertices of an $(n,d,\lambda)$-graph are removed independently with probability $1-p$ then, with high probability, the surviving graph is connected if $p\ge 1 - n^{-\,O(\frac{1}{d})}$ with the hidden constants depending only on $\frac{\lambda}{d}$.
\subsection{Related work}

In \citep{chung_giant_2009} the bond percolation model with a uniform
edge survival probability $p$ is studied. With $d_{i}$ denoting
the degree of vertex $i$, it is shown in \citep{chung_giant_2009}
that asymptotically almost surely a giant component survives (or not)
if $p>(1+\epsilon)\frac{\sum_{i}d_{i}}{\sum_{i}d_{i}^{2}}$ (or $p<(1-\epsilon)\frac{\sum_{i}d_{i}}{\sum_{i}d_{i}^{2}}$).
Furthermore, the spectral gap under bond percolation is studied in
\citep{chung_spectral_2007} and \citep{oliveira_concentration_2009},
where the latter established a sharper bound by means of concentration
inequalities for random matrices.

 A more relevant problem to our work
is the problem of\emph{ network (un)reliability} \citep{colbourn_combinatorics_1987}
where the goal is to estimate the probability that a percolated graph
remains connected. Under the bond percolation model, \citep{karger_randomized_2001} proposes a method to approximate the network reliability through a fully polynomial-time
approximation scheme. Approximation algorithms for the same problem with better computational
complexity were proposed later in \citep{harris_improved_2014} and \citep{karger_fast_2016}. 

The site percolation model for random $d$-regular graphs is analyzed
in \citep{greenhill_expansion_2008}. Specifically, \citep{greenhill_expansion_2008}
shows that, with high probability, for vertex deletion probability of
the form $n^{-\gamma}$, the surviving subgraph has a giant component
of order $n-o(n)$ that is an \emph{expander graph} and, if $\gamma\ge\frac{1}{d-1}$, then 
it is connected as well. This result was later improved and generalized
in \citep{ben-shimon_vertex_2009}. Recall that an $(n,d,\lambda)$-graph is a $d$-regular graph of order
$n$ with the non-trivial eigenvalue less than $\lambda$ in magnitude. A phase transition for site percolation on such $(n,d,\lambda)$-graphs 
is established in \citep{krivelevich_phase_2016}. In particular,
the mentioned paper shows that if the vertex survival probability
is $p=\frac{1-\epsilon}{d}$, then with high probability, the connected
components of the surviving subgraph have $O\left(\log n\right)$
vertices; whereas if $p=\frac{1+\epsilon}{d}$, $d=o(n)$, and $\frac{\lambda}{d}$
is relatively small, then with high probability a giant component
with $\Omega(\frac{n}{d})$ vertices survives.

Our main result, Theorem \ref{thm:aug-Laplacian-tail}, relies on a refined concentration inequality
stated in Proposition \ref{prp:Bernoulli-series-tail} for random
Bernoulli matrices and, consequently, is distinct from most of the previous
work mentioned above which rely on combinatorial arguments. We also
apply our general result to the special case of $(n,d,\lambda)$-graphs
(Corollary \ref{cor:ndl-graphs}), and reproduce bounds comparable
to those established in \citep{greenhill_expansion_2008,ben-shimon_vertex_2009}.
In particular, \citep[Proposition 3.5]{ben-shimon_vertex_2009} shows
that any $(n,d,\lambda)$-graph $G$ with $d\ge3$ and $\lambda\le2\sqrt{d-1}+\frac{1}{40}$,
that is also ``locally sparse'' in the sense that
\begin{align*}
\max_{\substack{H\subseteq G\,,\\
\left|V(H)\right|\le d+29
}
}\frac{\left|E(H)\right|}{\left|V(H)\right|} & \le1\,,
\end{align*}
 with high probability, remains connected under a homogeneous site
percolation with survival probability $p>1-n^{-\frac{1}{d}}$. Similarly, our result in Corollary \ref{cor:ndl-graphs}, guarantees that with probability $\ge1 - \frac{4}{n}$,
any $\left(n,d,\lambda\right)$-graph remains connected under a homogeneous
site percolation with survival probability $p\ge1-C_{1}n^{-\frac{2 C_{2}}{d}}$.
The constants $C_{1}$ and $C_{2}$ depend only  on  $\frac{\lambda}{d}$; their exact form is provided in the proof of Corollary \ref{cor:ndl-graphs}. If we have $\frac{\lambda}{d} = 1-\epsilon$ for some $\epsilon\in(0,1)$, then $-\log C_1= O\left((1+\frac{2}{\epsilon})^2\right)$ and $C_2 = O\left( \epsilon^{-4}\right)$ both of which are decreasing in $\epsilon$. These quantities can be fairly large for small values of $\epsilon$, which implies that our required lower bound on $p$ would be stricter than that of \citep{ben-shimon_vertex_2009}. However, our analysis does not explicitly
assume a bound on $\lambda$ or local sparsity as in \citep{ben-shimon_vertex_2009}.
The fact that Corollary \ref{cor:ndl-graphs} leads to suboptimal constants
compared to \citep[Proposition 3.5]{ben-shimon_vertex_2009} is not
surprising; Corollary \ref{cor:ndl-graphs} is an application of a very
general bound established in Theorem \ref{thm:aug-Laplacian-tail}
to the case of $(n,d,\lambda)$-graphs.

\subsection{Future directions}

There are natural extensions to the connectivity problem studied in
this paper that we would like to study through the lens of random
matrix theory as done here. For example, an immediate question is
to find a bound on the size of the giant component of the site-percolated
graph. Furthermore, an interesting research direction is to study
other properties of the site-percolated random graphs such as their
\emph{clique number}, \emph{chromatic number},\emph{ }etc by means
of random matrix theory. While the best results might still be obtained
through specifically tailored combinatorial arguments, we believe
that the analysis based on algebraic methods and random matrix theory
would be more robust to model errors.

\section{Problem Setup}

For $1\le i\le n$, let $\delta_{i}\sim\mr{Bernoulli}(p_{i})$ be
the independent random variables that indicate whether or not the
corresponding vertex survives. In order to operate on a Laplacian with fixed dimensions we interpret site percolation as removing every edge connected to the affected vertices. The Laplacian of the remaining graph $G_{\mb{\delta}}$ is then given by
\begin{align*}
\boldsymbol{L}_{\mb{\delta}} & =\sum_{1\le i<j\le n}\delta_{i}\delta_{j}w_{i,j}(\mb e_{i}-\mb e_{j})(\mb e_{i}-\mb e_{j})^{\T},
\end{align*}
 which also includes the vertices affected by the site percolation as isolated vertices. We need to take into account the effect of these ``ghost vertices'' to
find a non-trivial bound for the desired algebraic connectivity which
we denote by $a_{\mb{\delta}}$. To this end, for a coefficient $\alpha\ge0$, we introduce the
\emph{augmented Laplacian} given by 
\begin{align}
\overline{\mb L}_{\mb{\delta}} & =\mb L_{\mb{\delta}}+\sum_{1\le i\le n}\alpha\left(1-\delta_{i}\right)\mb e_{i}\mb e_{i}^{\T}.\label{eq:aug_L}
\end{align}
 The Laplacian $\mb L_{\mb{\delta}}$ and the diagonal matrix $\sum_{1\le i\le n}\alpha\left(1-\delta_{i}\right)\mb e_{i}\mb e_{i}^{\T}$ in \eqref{eq:aug_L} are supported on the vertices that survived and the ghost vertices, respectively. Because these two vertex sets are disjoint, the eigenvalues and
eigenvectors of the corresponding terms on the right-hand side of \eqref{eq:aug_L} constitute the eigendecomposition of $\overline{\mb L}_{\mb{\delta}}$.
We either have $a_{\mb{\delta}}>\alpha$ or $a_{\mb{\delta}}\le\alpha$.
If the latter holds, then $a_{\mb{\delta}}$ would coincide with the second smallest eigenvalue of $\overline{\mb{L}}_{\mb{\delta}}$ and by Weyl's eigenvalues inequality we obtain
\begin{align*}
a_{\mb{\delta}} & =\lambda_{2}\left(\overline{\mb L}_{\mb{\delta}}\right)\ge\lambda_{2}\left(\E\overline{\mb L}_{\mb{\delta}}\right)-\left\lVert \overline{\mb L}_{\mb{\delta}}-\E\overline{\mb L}_{\mb{\delta}}\right\rVert .
\end{align*}
 An immediate implication is that
\begin{align}
a_{\mb{\delta}} & \ge\min\left\{ \lambda_{2}\left(\E\overline{\mb L}_{\mb{\delta}}\right)-\left\lVert \overline{\mb L}_{\mb{\delta}}-\E\overline{\mb L}_{\mb{\delta}}\right\rVert ,\alpha\right\} ,\label{eq:a_delta_LB}
\end{align}
 holds for all $\alpha\ge0$. Hence, we can obtain a non-trivial lower
bound for $a_{\mb{\delta}}$ by studying the tail behavior of $\left\lVert \overline{\mb L}_{\mb{\delta}}-\E\overline{\mb L}_{\mb{\delta}}\right\rVert $
which also depends on $\alpha$. The lower bound given by (\ref{eq:a_delta_LB})
can also be optimized with respect to $\alpha$.

\section{Main Result}

Our main theorem below provides an upper bound for $\left\lVert \overline{\mb L}_{\mb{\delta}}-\E\overline{\mb L}_{\mb{\delta}}\right\rVert $.
To state the theorem it is necessary to introduce some notation. For
each $1\le i\le n$, let 
\begin{align}
K_{i} & =\frac{1}{2}\sqrt{\frac{1-2p_{i}}{\log\frac{1-p_{i}}{p_{i}}}}\label{eq:KS-var}
\end{align}
 denote the sub-Gaussian parameter of $\delta_{i}-p_{i}$ as used
in the Kearns-Saul inequality (Lemma \ref{lem:KS} below). Compared to the bounds on the moment generating function used in the Hoeffding and the Bernstein inequalities, the parameter \eqref{eq:KS-var} yields tighter bounds, particularly, if $p_i$ is close to $0$ or $1$. This property is crucial in our analysis as non-trivial events occur if the vertex survival probabilities (i.e., $p_i$s) are relatively close to $1$. We use
$\mb p=\left[\begin{array}{cccc}
 p_{1} & p_{2} & \dotsm & p_{n}\end{array}\right]^{\T}$ to denote the vector of survival probabilities and $\mb a_{i}$ to
denote the $i$th column of the adjacency matrix $\mb A$. The diagonal matrix whose diagonal entries are given by a vector $\mb{u}$ is denoted by $\mb{D}_{\mb{u}}$. The binary operation $\circ$ denotes the entrywise (or Hadamard) product.

\begin{thm}
\label{thm:aug-Laplacian-tail} Let $\delta_{i}\sim\mr{Bernoulli}(p_{i})$
be independent random variables for $1\le i\le n$. Furthermore, with
$K_{i}$ given by (\ref{eq:KS-var}) define 
\ifonecol
	\begin{align*}
	\sigma^{2} & =\left\lVert \sum_{i}K_{i}^{2}(1-2p_{i})^{2}\left(\mb a_{i}\mb a_{i}^{\T}\right)^{2}\right\rVert , & \overline{K} & =\max_{i}\left(\sum_{j}w_{i,j}^{2}K_{j}^{2}\right)^{\frac{1}{2}}.
	\end{align*}
\else
	\begin{align*}
	\sigma^{2} & =\left\lVert \sum_{i}K_{i}^{2}(1-2p_{i})^{2}\left(\mb a_{i}\mb a_{i}^{\T}\right)^{2}\right\rVert ,\\  \overline{K} & =\max_{i}\left(\sum_{j}w_{i,j}^{2}K_{j}^{2}\right)^{\frac{1}{2}}.
	\end{align*}
\fi

 Then for any $\varepsilon\in(0,1)$, with probability $\ge1-\varepsilon$
we have
\ifonecol
	\begin{align*}
	\left\lVert \overline{\mb L}_{\mb{\delta}}-\E\overline{\mb L}_{\mb{\delta}}\right\rVert  & \le2\overline{K}\sqrt{\log\frac{4n}{\varepsilon}}+\max_{i}\left|\alpha-\sum_{j}p_{j}w_{i,j}\right|\\
	 & \hspace{1em} + \left\lVert \mb D_{\mb p\circ\left(\mb 1-\mb p\right)}^{\frac{1}{2}}\mb A\mb D_{\mb p\circ\left(\mb 1-\mb p\right)}^{\frac{1}{2}}\right\rVert +2\left\lVert \mb D_{\mb p}\mb A\mb D_{\mb p\circ(\mb 1-\mb p)}^{\frac{1}{2}}\right\rVert +\frac{9}{2}\sqrt{\sigma\left(\log\frac{4n}{\varepsilon}\right)^{\frac{1}{2}}}.
	\end{align*}
\else
	\begin{align*}
	\left\lVert \overline{\mb L}_{\mb{\delta}}-\E\overline{\mb L}_{\mb{\delta}}\right\rVert  & \le2\overline{K}\sqrt{\log\frac{4n}{\varepsilon}}+\max_{i}\left|\alpha-\sum_{j}p_{j}w_{i,j}\right|\\
	 & \hspace{1em}+\left\lVert \mb D_{\mb p\circ\left(\mb 1-\mb p\right)}^{\frac{1}{2}}\mb A\mb D_{\mb p\circ\left(\mb 1-\mb p\right)}^{\frac{1}{2}}\right\rVert\\ & \hspace{1em} +2\left\lVert \mb D_{\mb p}\mb A\mb D_{\mb p\circ(\mb 1-\mb p)}^{\frac{1}{2}}\right\rVert +\frac{9}{2}\sqrt{\sigma\left(\log\frac{4n}{\varepsilon}\right)^{\frac{1}{2}}}.
	\end{align*}
\fi
\end{thm}
To evaluate the bound produced using (\ref{eq:a_delta_LB}) and Theorem
\ref{thm:aug-Laplacian-tail}, we apply the result to two special
problems with $\left(n,d,\lambda\right)$-graphs. First we recall the definition of these graphs.
\begin{defn}[$\left(n,d,\lambda\right)$-graphs]
 An $\left(n,d,\lambda\right)$-graph is a $d$-regular graph with
$n$ vertices whose adjacency matrix has no non-trivial eigenvalue
with magnitude greater than $\lambda$. 
\end{defn}
Below, we assume that the vertex deletion probabilities are identical,
i.e., $p_{1}=p_{2}=\dotsc=p_{n}=p$. This assumption also implies
that $K_{1}=K_{2}=\dotsc=K_{n}=K=\frac{1}{2}\sqrt{\frac{1-2p}{\log\frac{1-p}{p}}}$.
Also we assume all the edge weights $w_{i,j}$ are $\left\{ 0,1\right\} $-valued
and effectively indicate existence of an edge in $G$. We need to
quantify or bound $\lambda_{2}\left(\E\overline{\mb L}_{\mb{\delta}}\right)$
as well as the parameters $\sigma$ and $\overline{K}$. Using Theorem~\ref{thm:aug-Laplacian-tail}, the following corollary basically shows
that $p=1-\left(\frac{4n}{\varepsilon}\right)^{-O(\frac{1}{d})}$\,,
could suffice for any $(n,d,\lambda)$-graph affected by the prescribed
site percolation to remain connected with probability $\ge1-\varepsilon$. 
\begin{cor}
\label{cor:ndl-graphs}Let $G$ be an arbitrary $(n,d,\lambda)$-graph. There are positive constants $C_{1}$ and $C_{2}$ depending
only on $\frac{\lambda}{d}$ such that under the site percolation
model with vertex survival probability of 
\begin{align*}
p & \ge1-C_{1}\left(\frac{4n}{\varepsilon}\right)^{-C_{2}/d}
\end{align*}
 the surviving subgraph of $G$ is connected with probability $\ge1-\varepsilon$. 
\end{cor}
\begin{IEEEproof}
With $\mb A$ and $\mb L$ denoting the adjacency and Laplacian matrices
of the $\left(n,d,\lambda\right)$-graph $G$, the expected value
of the augmented Laplacian under the considered site percolation would
be
\begin{align*}
\E\overline{\mb L}_{\mb{\delta}} & =p^{2}\mb L+\alpha\left(1-p\right)\mb I=\left(p^{2}d+\alpha\left(1-p\right)\right)\mb I-p^{2}\mb A.
\end{align*}
 Let $\alpha=pd$. It follows from the definition of the graph and
the equation above that
\begin{align}
\lambda_{2}\left(\E\overline{\mb L}_{\mb{\delta}}\right) & \ge p^{2}\left(d-\lambda\right)+p(1-p)d=pd-p^{2}\lambda\,.\label{eq:Elambda_ndl}
\end{align}
 Furthermore, the parameters $\sigma^{2}$ and $\overline{K}$ can
be expressed as
\ifonecol
	\begin{align*}
	\sigma^{2} & =K^{2}\left(1-2p\right)^{2}\left\lVert \sum_{i}\left(\mb a_{i}\mb a_{i}^{\T}\right)^{2}\right\rVert  & \text{and} &  & \overline{K} & =K\sqrt{d}.\\
	 & =K^{2}\left(1-2p\right)^{2}d\left\lVert \mb A^{2}\right\rVert =K^{2}\left(1-2p\right)^{2}d^{3}
	\end{align*}
\else
	\begin{align*}
	\sigma^{2} & =K^{2}\left(1-2p\right)^{2}\left\lVert \sum_{i}\left(\mb a_{i}\mb a_{i}^{\T}\right)^{2}\right\rVert \\
	 & =K^{2}\left(1-2p\right)^{2}d\left\lVert \mb A^{2}\right\rVert =K^{2}\left(1-2p\right)^{2}d^{3}\\
	\end{align*}
	and 
	\begin{align*}
	 \overline{K} & =K\sqrt{d}.
	\end{align*}
\fi
 Finally, we have $\left\lVert \mb D_{\mb p\circ\left(\mb 1-\mb p\right)}^{\frac{1}{2}}\mb A\mb D_{\mb p\circ\left(\mb 1-\mb p\right)}^{\frac{1}{2}}\right\rVert =p\left(1-p\right)\left\lVert \mb A\right\rVert =p\left(1-p\right)d$,
$\left\lVert \mb D_{\mb p}\mb A\mb D_{\mb p\circ(\mb 1-\mb p)}^{\frac{1}{2}}\right\rVert =p^{\frac{3}{2}}\left(1-p\right)^{\frac{1}{2}}d$,
and 
\begin{align*}
\max_{i}\left|\alpha-\sum_{j}p_{j}w_{i,j}\right| & =0.
\end{align*}

We can now invoke Theorem \ref{thm:aug-Laplacian-tail} and apply
the above bounds to obtain 
\ifonecol
	\begin{align*}
	\left\lVert \overline{\mb L}_{\mb{\delta}}-\E\overline{\mb L}_{\mb{\delta}}\right\rVert  & \le2K\sqrt{d\log\frac{4n}{\varepsilon}}+p\left(1-p\right)d+2p^{\frac{3}{2}}\left(1-p\right)^{\frac{1}{2}}d+\frac{9}{2}\sqrt{K\left|1-2p\right|}d^{\frac{3}{4}}\left(\log\frac{4n}{\varepsilon}\right)^{\frac{1}{4}}.
	\end{align*}
\else
	\begin{align*}
\left\lVert \overline{\mb L}_{\mb{\delta}}-\E\overline{\mb L}_{\mb{\delta}}\right\rVert  & \le2K\sqrt{d\log\frac{4n}{\varepsilon}}+p\left(1-p\right)d+2p^{\frac{3}{2}}\left(1-p\right)^{\frac{1}{2}}d\\
&\hspace{1em}+\frac{9}{2}\sqrt{K\left|1-2p\right|}d^{\frac{3}{4}}\left(\log\frac{4n}{\varepsilon}\right)^{\frac{1}{4}}.
\end{align*}
\fi
 Given the inequalities (\ref{eq:a_delta_LB}), (\ref{eq:Elambda_ndl}),
and the assumption that $\alpha=pd$, we are naturally interested
in values of $p$ for which the right-hand side of the inequality
above is strictly smaller than $pd-p^{2}\lambda$. Specifically, we
would like to find $p$ for which we have
\ifonecol
	\begin{align*}
	pd-p^{2}\lambda & >2K\sqrt{d\log\frac{4n}{\varepsilon}}+p\left(1-p\right)d+2p^{\frac{3}{2}}\left(1-p\right)^{\frac{1}{2}}d+\frac{9}{2}\sqrt{K\left|1-2p\right|}d^{\frac{3}{4}}\left(\log\frac{4n}{\varepsilon}\right)^{\frac{1}{4}}\,,
	\end{align*}
 	or equivalently
	\begin{align}
	(1-\frac{\lambda}{d})p & >2\left(\frac{\log\frac{4n}{\varepsilon}}{d}\cdot K^{2}\right)^{\frac{1}{2}}+2\sqrt{p(1-p)}+\frac{9}{2p}\sqrt{\left|1-2p\right|}\left(\frac{\log\frac{4n}{\varepsilon}}{d}\cdot K^{2}\right)^{\frac{1}{4}}.\label{eq:ndlambda-main}
	\end{align}
\else
	\begin{align*}
	pd-p^{2}\lambda & >2K\sqrt{d\log\frac{4n}{\varepsilon}}+p\left(1-p\right)d+2p^{\frac{3}{2}}\left(1-p\right)^{\frac{1}{2}}d\\
	&\hspace{1em}+\frac{9}{2}\sqrt{K\left|1-2p\right|}d^{\frac{3}{4}}\left(\log\frac{4n}{\varepsilon}\right)^{\frac{1}{4}}\,,
	\end{align*}
	 or equivalently
	\begin{align}
	\begin{aligned}
	(1-\frac{\lambda}{d})p & >2\left(\frac{\log\frac{4n}{\varepsilon}}{d}\cdot K^{2}\right)^{\frac{1}{2}}+2\sqrt{p(1-p)}\\
	& \hspace{1em}+\frac{9}{2p}\sqrt{\left|1-2p\right|}\left(\frac{\log\frac{4n}{\varepsilon}}{d}\cdot K^{2}\right)^{\frac{1}{4}}.
	\end{aligned}
	\label{eq:ndlambda-main}
	\end{align}
\fi
 For $p\ge\frac{1}{2}$ we have $K\le\frac{1}{2\sqrt{\log\frac{1}{1-p}}}$.
Therefore, if we parameterize $p$ by $\beta\ge1$ as $p=1-e^{-\beta^{4}}$
we have $K^{2}\le\frac{\beta^{-4}}{4}$. Furthermore, we can write
$p\ge\frac{1}{1+\beta^{-4}}\ge1-\beta^{-2}$, $2\sqrt{p(1-p)}\le2e^{-\beta^{4}/2}\le2\beta^{-2}$,
and $\frac{1}{2p}\sqrt{\left|1-2p\right|}\le1$. Therefore, to guarantee
(\ref{eq:ndlambda-main}) it suffices to have 
\ifonecol
	\begin{align*}
	(1-\frac{\lambda}{d})(1-\beta^{-2}) & >2\left(\frac{\beta^{-4}\log\frac{4n}{\varepsilon}}{4d}\right)^{\frac{1}{2}}+2\beta^{-2}+9\left(\frac{\beta^{-4}\log\frac{4n}{\varepsilon}}{4d}\right)^{\frac{1}{4}}\\
	 & =\left(\left(\frac{\log\frac{4n}{\varepsilon}}{d}\right)^{\frac{1}{2}}+2\right)\beta^{-2}+\frac{9}{\sqrt{2}}\left(\frac{\log\frac{4n}{\varepsilon}}{d}\right)^{\frac{1}{4}}\beta^{-1}\,,
	\end{align*}
\else
	\begin{align*}
	(1\!-\!\frac{\lambda}{d})(1\!-\!\beta^{-2}) & >2\left(\!\frac{\beta^{-4}\log\frac{4n}{\varepsilon}}{4d}\right)^{\hspace{-3pt}\frac{1}{2}}+2\beta^{-2}+9\left(\!\frac{\beta^{-4}\log\frac{4n}{\varepsilon}}{4d}\right)^{\hspace{-3pt}\frac{1}{4}}\\
	 & =\left(\!\left(\frac{\log\frac{4n}{\varepsilon}}{d}\right)^{\hspace{-3pt}\frac{1}{2}}\hspace{-5pt}+2\right)\beta^{-2}+\frac{9}{\sqrt{2}}\left(\frac{\log\frac{4n}{\varepsilon}}{d}\right)^{\hspace{-3pt}\frac{1}{4}}\hspace{-5pt}\beta^{-1}\,,
	\end{align*}
\fi
 which is equivalent to 
\begin{align*}
1-\frac{\lambda}{d} & >\left(\left(\frac{\log\frac{4n}{\varepsilon}}{d}\right)^{\frac{1}{2}}+3-\frac{\lambda}{d}\right)\beta^{-2}+\frac{9}{\sqrt{2}}\left(\frac{\log\frac{4n}{\varepsilon}}{d}\right)^{\frac{1}{4}}\beta^{-1}.
\end{align*}
The inequality above holds for
\ifonecol
	\begin{align*}
	\beta^{2} & >\max\left\{ 2\left(1-\frac{\lambda}{d}\right)^{-1}\left(\left(\frac{\log\frac{4n}{\varepsilon}}{d}\right)^{\frac{1}{2}}+3-\frac{\lambda}{d}\right),81\left(1-\frac{\lambda}{d}\right)^{-2}\left(\frac{\log\frac{4n}{\varepsilon}}{d}\right)^{\frac{1}{2}}\right\} 
	\end{align*}
\else
	\begin{align*}
	\beta^{2} & >\max\left\lbrace 2\left(1-\frac{\lambda}{d}\right)^{-1}\left(\left(\frac{\log\frac{4n}{\varepsilon}}{d}\right)^{\frac{1}{2}}+3-\frac{\lambda}{d}\right) , \right. \\ & \hspace{5em}\left. 81\left(1-\frac{\lambda}{d}\right)^{-2}\left(\frac{\log\frac{4n}{\varepsilon}}{d}\right)^{\frac{1}{2}}\right\rbrace 
	\end{align*}
\fi
 which, for the sake of simpler expressions, can be further relaxed
to 
\begin{align*}
\beta^{4} & \ge 81^{2}\left(1-\frac{\lambda}{d}\right)^{-4}\frac{\log\frac{4n}{\varepsilon}}{d}+8\left(1-\frac{\lambda}{d}\right)^{-2}\left(3-\frac{\lambda}{d}\right)^{2}.
\end{align*}
 The desired results follows immediately by setting $C_{1}=\exp\left(-8\left(1-\frac{\lambda}{d}\right)^{-2}\left(3-\frac{\lambda}{d}\right)^{2}\right)$
and $C_{2}=81^{2}\left(1-\frac{\lambda}{d}\right)^{-4}$.
\end{IEEEproof}

\section{Proof of Theorem \ref{thm:aug-Laplacian-tail}}

In this section we prove our main result. The lemmas and other technical
tools we use are summarized below in Appendix \ref{apx:lemmas}.
\begin{IEEEproof}[\textbf{Proof of Theorem }\ref{thm:aug-Laplacian-tail}]
Splitting $\overline{\mb L}_{\mb{\delta}}-\E\overline{\mb L}_{\mb{\delta}}$
into the sum of diagonal and off-diagonal terms as 
\ifonecol
	\begin{align*}
	\overline{\mb L}_{\mb{\delta}}-\E\overline{\mb L}_{\mb{\delta}} & =\sum_{i}\left(\delta_{i}(\sum_{j}\delta_{j}w_{i,j})-p_{i}(\sum_{j}p_{j}w_{i,j})-\alpha(\delta_{i}-p_{i})\right)\mb e_{i}\mb e_{i}^{\T}\\
	 & \hspace{1em}+\sum_{i<j}(\delta_{i}\delta_{j}-p_{i}p_{j})(\mb e_{i}\mb e_{j}^{\T}+\mb e_{j}\mb e_{i}^{\T})
	\end{align*}
\else
	\begin{align*}
	\overline{\mb L}_{\mb{\delta}}\!-\!\E\overline{\mb L}_{\mb{\delta}} & \!=\!\sum_{i}\!\left(\!\delta_{i}(\sum_{j}\!\delta_{j}w_{i,j})-p_{i}(\sum_{j}p_{j}w_{i,j})-\alpha(\delta_{i}\!-\! p_{i})\!\!\right)\!\mb e_{i}\mb e_{i}^{\T}\\
	 & \hspace{1em}+\sum_{i<j}(\delta_{i}\delta_{j}-p_{i}p_{j})(\mb e_{i}\mb e_{j}^{\T}+\mb e_{j}\mb e_{i}^{\T})
	\end{align*}
\fi
and applying triangle inequality yields
\ifonecol
	\begin{align*}
	\left\lVert \overline{\mb L}_{\mb{\delta}}-\E\overline{\mb L}_{\mb{\delta}}\right\rVert  & \le\max_{i}\left|\delta_{i}(\sum_{j}\delta_{j}w_{i,j})-p_{i}(\sum_{j}p_{j}w_{i,j})-\alpha(\delta_{i}-p_{i})\right|\\
	 & \hspace{1em}+\left\lVert \sum_{i<j}(\delta_{i}\delta_{j}-p_{i}p_{j})w_{i,j}(\mb e_{i}\mb e_{j}^{\T}+\mb e_{j}\mb e_{i}^{\T})\right\rVert .
	\end{align*}
\else
	\begin{align*}
	\left\lVert \overline{\mb L}_{\mb{\delta}}-\E\overline{\mb L}_{\mb{\delta}}\right\rVert  & \le\max_{i}\left|\delta_{i}(\sum_{j}\!\delta_{j}w_{i,j})-p_{i}(\sum_{j}p_{j}w_{i,j})-\alpha(\delta_{i}\!-\! p_{i})\right|\\
	 & \hspace{1em}+\left\lVert \sum_{i<j}(\delta_{i}\delta_{j}-p_{i}p_{j})w_{i,j}(\mb e_{i}\mb e_{j}^{\T}+\mb e_{j}\mb e_{i}^{\T})\right\rVert .
	\end{align*}
\fi
Our goal is to bound the two terms on the right-hand side of the inequality
above. To lighten the notation we use $\xi_{i}=\delta_{i}-p_{i}$
for $i=1,2,\dotsc,n$ and let
\begin{align*}
S_{1} & =\max_{i}\left|\delta_{i}(\sum_{j}\delta_{j}w_{i,j})-p_{i}(\sum_{j}p_{j}w_{i,j})-\alpha(\delta_{i}-p_{i})\right|
\end{align*}
 and 
\begin{align*}
S_{2} & =\left\lVert \sum_{i<j}(\delta_{i}\delta_{j}-p_{i}p_{j})w_{i,j}(\mb e_{i}\mb e_{j}^{\T}+\mb e_{j}\mb e_{i}^{\T})\right\rVert .
\end{align*}

It is straightforward to verify that
\ifonecol
	\begin{align*}
	\left|\delta_{i}(\sum_{j}\delta_{j}w_{i,j})-p_{i}(\sum_{j}p_{j}w_{i,j})-\alpha(\delta_{i}-p_{i})\right| & =\left|\delta_{i}\sum_{j}(\delta_{j}-p_{j})w_{i,j}-(\delta_{i}-p_{i})(\alpha-\sum_{j}p_{j}w_{i,j})\right|\\
	 & \le\left|\sum_{j}\xi_{j}w_{i,j}\right|+\left|\alpha-\sum_{j}p_{j}w_{i,j}\right|
	\end{align*}
\else
	\begin{align*}
	& \left|\delta_{i}(\sum_{j}\delta_{j}w_{i,j})-p_{i}(\sum_{j}p_{j}w_{i,j})-\alpha(\delta_{i}-p_{i})\right|\\  = & \left|\delta_{i}\sum_{j}(\delta_{j}-p_{j})w_{i,j}-(\delta_{i}-p_{i})(\alpha-\sum_{j}p_{j}w_{i,j})\right|\\
	  \le & \left|\sum_{j}\xi_{j}w_{i,j}\right|+\left|\alpha-\sum_{j}p_{j}w_{i,j}\right|
	\end{align*}
\fi
 using which we obtain 
\begin{align*}
S_{1} & \le\max_{i}\left|\sum_{j}\xi_{j}w_{i,j}\right|+\left|\alpha-\sum_{j}p_{j}w_{i,j}\right|.
\end{align*}
 By Chernoff's inequality and Lemma \ref{lem:KS}, for each $i$ we
have 
\begin{align*}
\left|\sum_{j}\xi_{j}w_{i,j}\right| & \le2\left(\sum_{j}w_{i,j}^{2}K_{j}^{2}\right)^{\frac{1}{2}}\sqrt{\log\frac{4n}{\varepsilon}},
\end{align*}
with probability $\ge1-\frac{\varepsilon}{2n}$. Then by union bound
we have 
\ifonecol
	\begin{align}
	S_{1} & \le\max_{i}\left(\sum_{j}w_{i,j}^{2}K_{j}^{2}\right)^{\frac{1}{2}}\sqrt{\log\frac{8n}{\varepsilon}}+\left|\alpha-\sum_{j}p_{j}w_{i,j}\right|\le2\overline{K}\sqrt{\log\frac{4n}{\varepsilon}}+\max_{i}\left|\alpha-\sum_{j}p_{j}w_{i,j}\right|,\label{eq:S1<}
	\end{align}
\else
	\begin{align}
	S_{1} & \le\max_{i}\left(\sum_{j}w_{i,j}^{2}K_{j}^{2}\right)^{\frac{1}{2}}\sqrt{\log\frac{8n}{\varepsilon}}+\left|\alpha-\sum_{j}p_{j}w_{i,j}\right|\nonumber\\
	& \le2\overline{K}\sqrt{\log\frac{4n}{\varepsilon}}+\max_{i}\left|\alpha-\sum_{j}p_{j}w_{i,j}\right|,\label{eq:S1<}
	\end{align}
\fi
 with probability $\ge1-\frac{\varepsilon}{2}$.

Expressing $S_{2}$ in terms of $\xi_{i}$s and applying the triangle
inequality reveals that
\ifonecol 
	\begin{align}
	S_{2} & =\left\lVert \sum_{i<j}(\xi_{i}\xi_{j}-p_{i}\xi_{j}-\xi_{i}p_{j})w_{i,j}(\mb e_{i}\mb e_{j}^{\T}+\mb e_{j}\mb e_{i}^{\T})\right\rVert \nonumber \\
	 & \le\left\lVert \sum_{i<j}\xi_{i}\xi_{j}w_{i,j}(\mb e_{i}\mb e_{j}^{\T}+\mb e_{j}\mb e_{i}^{\T})\right\rVert +\left\lVert \sum_{i<j}\left(p_{i}\xi_{j}+p_{j}\xi_{i}\right)w_{i,j}(\mb e_{i}\mb e_{j}^{\T}+\mb e_{j}\mb e_{i}^{\T})\right\rVert \nonumber \\
	 & =\left\lVert \mb D_{\mb{\xi}}\mb A\mb D_{\mb{\xi}}\right\rVert +2\left\lVert \mb D_{\mb p}\mb A\mb D_{\mb{\xi}}\right\rVert ,\label{eq:S2-bound}
	\end{align}
\else
	\begin{align}
	S_{2} & =\left\lVert \sum_{i<j}(\xi_{i}\xi_{j}-p_{i}\xi_{j}-\xi_{i}p_{j})w_{i,j}(\mb e_{i}\mb e_{j}^{\T}+\mb e_{j}\mb e_{i}^{\T})\right\rVert \nonumber \\
	 & \le\left\lVert \sum_{i<j}\xi_{i}\xi_{j}w_{i,j}(\mb e_{i}\mb e_{j}^{\T}+\mb e_{j}\mb e_{i}^{\T})\right\rVert\nonumber \\
	 & \hspace{1em} +\left\lVert \sum_{i<j}\left(p_{i}\xi_{j}+p_{j}\xi_{i}\right)w_{i,j}(\mb e_{i}\mb e_{j}^{\T}+\mb e_{j}\mb e_{i}^{\T})\right\rVert \nonumber \\
	 & =\left\lVert \mb D_{\mb{\xi}}\mb A\mb D_{\mb{\xi}}\right\rVert +2\left\lVert \mb D_{\mb p}\mb A\mb D_{\mb{\xi}}\right\rVert ,\label{eq:S2-bound}
	\end{align}
\fi
with $\mb D_{\mb{\xi}}$ and $\mb D_{\mb p}$ respectively denoting
diagonal matrices with $\xi_{i}$s and $p_{i}$s on their diagonals.

We can write
\ifonecol 
	\begin{align}
	\left\lVert \mb D_{\mb{\xi}}\mb A\mb D_{\mb{\xi}}\right\rVert ^{2} & =\left\lVert \mb D_{\mb{\xi}}\mb A\mb D_{\mb{\xi}}^{2}\mb A\mb D_{\mb{\xi}}\right\rVert \nonumber \\
	 & \le\left\lVert \mb D_{\mb{\xi}}\mb A\mb D_{\mb p\circ\left(\mb 1-\mb p\right)}\mb A\mb D_{\mb{\xi}}\right\rVert +\left\lVert \mb D_{\mb{\xi}}\mb A\mb D_{\mb{\xi}\circ(\mb 1-2\mb p)}\mb A\mb D_{\mb{\xi}}\right\rVert \nonumber \\
	 & \le\left\lVert \mb D_{\mb p\circ\left(\mb 1-\mb p\right)}^{\frac{1}{2}}\mb A\mb D_{\mb{\xi}}^{2}\mb A\mb D_{\mb p\circ\left(\mb 1-\mb p\right)}^{\frac{1}{2}}\right\rVert +\left\lVert \mb A\mb D_{\mb{\xi}\circ(\mb 1-2\mb p)}\mb A\right\rVert \nonumber \\
	 & \le\left\lVert \mb D_{\mb p\circ\left(\mb 1-\mb p\right)}^{\frac{1}{2}}\mb A\mb D_{\mb p\circ\left(\mb 1-\mb p\right)}\mb A\mb D_{\mb p\circ\left(\mb 1-\mb p\right)}^{\frac{1}{2}}\right\rVert \nonumber \\
	 & \hspace{1em}+\left\lVert \mb D_{\mb p\circ\left(\mb 1-\mb p\right)}^{\frac{1}{2}}\mb A\mb D_{\mb{\xi}\circ(\mb 1-2\mb p)}\mb A\mb D_{\mb p\circ\left(\mb 1-\mb p\right)}^{\frac{1}{2}}\right\rVert +\left\lVert \mb A\mb D_{\mb{\xi}\circ(\mb 1-2\mb p)}\mb A\right\rVert .\label{eq:DAD^2}
	\end{align}
\else
	\begin{align}
	\left\lVert \mb D_{\mb{\xi}}\mb A\mb D_{\mb{\xi}}\right\rVert ^{2} & =\left\lVert \mb D_{\mb{\xi}}\mb A\mb D_{\mb{\xi}}^{2}\mb A\mb D_{\mb{\xi}}\right\rVert \nonumber \\
	 & \le\left\lVert \mb D_{\mb{\xi}}\mb A\mb D_{\mb p\circ\left(\mb 1-\mb p\right)}\mb A\mb D_{\mb{\xi}}\right\rVert \nonumber \\
	 & \hspace{1em}+\left\lVert \mb D_{\mb{\xi}}\mb A\mb D_{\mb{\xi}\circ(\mb 1-2\mb p)}\mb A\mb D_{\mb{\xi}}\right\rVert \nonumber \\
	 & \le\left\lVert \mb D_{\mb p\circ\left(\mb 1-\mb p\right)}^{\frac{1}{2}}\mb A\mb D_{\mb{\xi}}^{2}\mb A\mb D_{\mb p\circ\left(\mb 1-\mb p\right)}^{\frac{1}{2}}\right\rVert \nonumber \\ & \hspace{1em}+\left\lVert \mb A\mb D_{\mb{\xi}\circ(\mb 1-2\mb p)}\mb A\right\rVert \nonumber \\
	 & \begin{aligned} 
	 & \le\left\lVert \mb D_{\mb p\circ\left(\mb 1-\mb p\right)}^{\frac{1}{2}}\mb A\mb D_{\mb p\circ\left(\mb 1-\mb p\right)}\mb A\mb D_{\mb p\circ\left(\mb 1-\mb p\right)}^{\frac{1}{2}}\right\rVert \\
	 & \hspace{1em}+\left\lVert \mb D_{\mb p\circ\left(\mb 1-\mb p\right)}^{\frac{1}{2}}\mb A\mb D_{\mb{\xi}\circ(\mb 1-2\mb p)}\mb A\mb D_{\mb p\circ\left(\mb 1-\mb p\right)}^{\frac{1}{2}}\right\rVert \\
	 & \hspace{1em}+\left\lVert \mb A\mb D_{\mb{\xi}\circ(\mb 1-2\mb p)}\mb A\right\rVert.
	 \end{aligned}
	  \label{eq:DAD^2}
	\end{align}
\fi
 We used the identity $\delta_{i}^{2}=\delta_{i}$ or equivalently
$\xi_{i}^{2}=p_{i}(1-p_{i})+\xi_{i}(1-2p_{i})$ followed by a triangle
inequality to obtain the first inequality.  To obtain the second inequality
we simply rearranged the matrices in the first term and used the fact
that $\left|\xi_{i}\right|\le1$ to bound the second term. Applying
the identity $\delta_{i}^{2}=\delta_{i}$ again yields the third inequality.
With $\mb a_{i}$ denoting the $i$th column of the adjacency matrix
$\mb A$, we can invoke Proposition \ref{prp:Bernoulli-series-tail}
to guarantee that with probability $\ge1-\frac{\varepsilon}{2}$ we
have 
\begin{align*}
\left\lVert \mb A\mb D_{\mb{\xi}\circ(\mb 1-2\mb p)}\mb A\right\rVert  & =\left\lVert \sum_{i}\xi_{i}\left(1-2p_{i}\right)\mb a_{i}\mb a_{i}^{\T}\right\rVert \\
 & \le2\left\lVert \sum_{i}K_{i}^{2}(1-2p_{i})^{2}\left(\mb a_{i}\mb a_{i}^{\T}\right)^{2}\right\rVert ^{\frac{1}{2}}\ifonecol \else \hspace{-5pt} \fi \sqrt{\log\frac{4n}{\varepsilon}}\\
 & =2\sigma\sqrt{\log\frac{4n}{\varepsilon}}.
\end{align*}
On the same event we also have
\ifonecol 
	\begin{align*}
	\left\lVert \mb D_{\mb p\circ\left(\mb 1-\mb p\right)}^{\frac{1}{2}}\mb A\mb D_{\mb{\xi}\circ(\mb 1-2\mb p)}\mb A\mb D_{\mb p\circ\left(\mb 1-\mb p\right)}^{\frac{1}{2}}\right\rVert  & \le\left\lVert \mb D_{\mb p\circ\left(\mb 1-\mb p\right)}^{\frac{1}{2}}\right\rVert \left\lVert \mb A\mb D_{\mb{\xi}\circ(\mb 1-2\mb p)}\mb A\right\rVert \left\lVert \mb D_{\mb p\circ\left(\mb 1-\mb p\right)}^{\frac{1}{2}}\right\rVert \\
	 & \le\frac{1}{4}\left\lVert \mb A\mb D_{\mb{\xi}\circ(\mb 1-2\mb p)}\mb A\right\rVert \\
	 & \le\frac{\sigma}{2}\sqrt{\log\frac{4n}{\varepsilon}}\,,
	\end{align*}
\else
	\begin{align*}
	& \left\lVert \mb D_{\mb p\circ\left(\mb 1-\mb p\right)}^{\frac{1}{2}}\mb A\mb D_{\mb{\xi}\circ(\mb 1-2\mb p)}\mb A\mb D_{\mb p\circ\left(\mb 1-\mb p\right)}^{\frac{1}{2}}\right\rVert \\ \le & \left\lVert \mb D_{\mb p\circ\left(\mb 1-\mb p\right)}^{\frac{1}{2}}\right\rVert \left\lVert \mb A\mb D_{\mb{\xi}\circ(\mb 1-2\mb p)}\mb A\right\rVert \left\lVert \mb D_{\mb p\circ\left(\mb 1-\mb p\right)}^{\frac{1}{2}}\right\rVert \\
	 \le & \frac{1}{4}\left\lVert \mb A\mb D_{\mb{\xi}\circ(\mb 1-2\mb p)}\mb A\right\rVert \\
	 \le & \frac{\sigma}{2}\sqrt{\log\frac{4n}{\varepsilon}}\,,
	\end{align*}
\fi
where we used the fact that $\left\lVert \mb D_{\mb p\circ\left(\mb 1-\mb p\right)}\right\rVert \le\frac{1}{4}$
in the second line. These probabilistic upper bounds together with
(\ref{eq:DAD^2}) ensure that\looseness=-1
\begin{align}
\left\lVert \mb D_{\mb{\xi}}\mb A\mb D_{\mb{\xi}}\right\rVert  & \le\left(\left\lVert \mb D_{\mb p\circ\left(\mb 1-\mb p\right)}^{\frac{1}{2}}\mb A\mb D_{\mb p\circ\left(\mb 1-\mb p\right)}^{\frac{1}{2}}\right\rVert ^{2}+\frac{5}{2}\sigma\sqrt{\log\frac{4n}{\varepsilon}}\right)^{\frac{1}{2}}\nonumber \\
 & \le\left\lVert \mb D_{\mb p\circ\left(\mb 1-\mb p\right)}^{\frac{1}{2}}\mb A\mb D_{\mb p\circ\left(\mb 1-\mb p\right)}^{\frac{1}{2}}\right\rVert +\sqrt{\frac{5}{2}\sigma\left(\log\frac{4n}{\varepsilon}\right)^{\frac{1}{2}}}\,,
\label{eq:DxiADxi}
\end{align}
with probability $\ge1-\frac{\varepsilon}{2}$. Furthermore, using
similar arguments as above we have
\ifonecol
	\begin{align*}
	\left\lVert \mb D_{\mb p}\mb A\mb D_{\mb{\xi}}\right\rVert ^{2} & =\left\lVert \mb D_{\mb p}\mb A\mb D_{\mb{\xi}}^{2}\mb A\mb D_{\mb p}\right\rVert \\
	 & \le\left\lVert \mb D_{\mb p}\mb A\mb D_{\mb p\circ(\mb 1-\mb p)}\mb A\mb D_{\mb p}\right\rVert +\left\lVert \mb D_{\mb p}\mb A\mb D_{\mb{\xi}\circ(\mb 1-2\mb p)}\mb A\mb D_{\mb p}\right\rVert \\
	 & \le\left\lVert \mb D_{\mb p}\mb A\mb D_{\mb p\circ(\mb 1-\mb p)}^{\frac{1}{2}}\right\rVert ^{2}+\left\lVert \mb D_{\mb p}\right\rVert ^{2}\left\lVert \mb A\mb D_{\mb{\xi}\circ(\mb 1-2\mb p)}\mb A\right\rVert \\
	 & \le\left\lVert \mb D_{\mb p}\mb A\mb D_{\mb p\circ(\mb 1-\mb p)}^{\frac{1}{2}}\right\rVert ^{2}+2\sigma\sqrt{\log\frac{4n}{\varepsilon}},
	\end{align*}
\else
	\begin{align*}
	\left\lVert \mb D_{\mb p}\mb A\mb D_{\mb{\xi}}\right\rVert ^{2} & =\left\lVert \mb D_{\mb p}\mb A\mb D_{\mb{\xi}}^{2}\mb A\mb D_{\mb p}\right\rVert \\
	 & \le\left\lVert \mb D_{\mb p}\mb A\mb D_{\mb p\circ(\mb 1-\mb p)}\mb A\mb D_{\mb p}\right\rVert \\
	 & \hspace{1em}+\left\lVert \mb D_{\mb p}\mb A\mb D_{\mb{\xi}\circ(\mb 1-2\mb p)}\mb A\mb D_{\mb p}\right\rVert \\
	 & \le\left\lVert \mb D_{\mb p}\mb A\mb D_{\mb p\circ(\mb 1-\mb p)}^{\frac{1}{2}}\right\rVert ^{2}\!\!+\!\left\lVert \mb D_{\mb p}\right\rVert ^{2}\left\lVert \mb A\mb D_{\mb{\xi}\circ(\mb 1-2\mb p)}\mb A\right\rVert \\
	 & \le\left\lVert \mb D_{\mb p}\mb A\mb D_{\mb p\circ(\mb 1-\mb p)}^{\frac{1}{2}}\right\rVert ^{2}+2\sigma\sqrt{\log\frac{4n}{\varepsilon}},
	\end{align*}
\fi
and thus 
\begin{align}
\left\lVert \mb D_{\mb p}\mb A\mb D_{\mb{\xi}}\right\rVert  & \le\left\lVert \mb D_{\mb p}\mb A\mb D_{\mb p\circ(\mb 1-\mb p)}^{\frac{1}{2}}\right\rVert +\sqrt{2\sigma\left(\log\frac{4n}{\varepsilon}\right)^{\frac{1}{2}}}.\label{eq:DpADxi}
\end{align}
 It follows from (\ref{eq:S2-bound}), (\ref{eq:DxiADxi}), and (\ref{eq:DpADxi})
that
\ifonecol
	\begin{align}
	S_{2} & \le\left\lVert \mb D_{\mb p\circ\left(\mb 1-\mb p\right)}^{\frac{1}{2}}\mb A\mb D_{\mb p\circ\left(\mb 1-\mb p\right)}^{\frac{1}{2}}\right\rVert +2\left\lVert \mb D_{\mb p}\mb A\mb D_{\mb p\circ(\mb 1-\mb p)}^{\frac{1}{2}}\right\rVert +\frac{9}{2}\sqrt{\sigma\left(\log\frac{4n}{\varepsilon}\right)^{\frac{1}{2}}},\label{eq:S2<}
	\end{align}
\else
	\begin{align}
	\begin{aligned}
	S_{2} & \le\left\lVert \mb D_{\mb p\circ\left(\mb 1-\mb p\right)}^{\frac{1}{2}}\mb A\mb D_{\mb p\circ\left(\mb 1-\mb p\right)}^{\frac{1}{2}}\right\rVert +2\left\lVert \mb D_{\mb p}\mb A\mb D_{\mb p\circ(\mb 1-\mb p)}^{\frac{1}{2}}\right\rVert \\
	& \hspace{1em}+\frac{9}{2}\sqrt{\sigma\left(\log\frac{4n}{\varepsilon}\right)^{\frac{1}{2}}},
	\end{aligned}\label{eq:S2<}
	\end{align}
\fi
with probability $\ge1-\frac{\varepsilon}{2}$. 

The desired result follows immediately using the derived bounds (\ref{eq:S1<})
and (\ref{eq:S2<}).
\end{IEEEproof}
\appendices

\section{\label{apx:lemmas}Auxiliary tools and technical lemmas}

We use the following lemma due to \citep{kearns_large_98} which provides
a sharp bound for the sub-Gaussian norm of general Bernoulli random
variables.
\begin{lem}[Kearns-Saul inequality \citep{kearns_large_98}]
\label{lem:KS} For $p\in\left[0,1\right]$ let $\delta$ be a $\mr{Bernoulli}(p)$
random variable. Then for all $t\in\mbb R$ we have 
\begin{align*}
\E e^{t\left(\delta-p\right)}=pe^{t(1-p)}+(1-p)e^{-tp} & \le e^{\left(K(p)\,t\right)^{2}},
\end{align*}
where 
\begin{align*}
K(p) & \defeq\frac{1}{2}\sqrt{\frac{1-2p}{\log\frac{1-p}{p}}}.
\end{align*}
\end{lem}
We also use the following master tail bound for sums of independent
random matrices due to \citep{tropp_user-friendly_2011}. 
\begin{thm}
\label{thm:master-tail}\citep[Theorem 3.6]{tropp_user-friendly_2011}
Consider a finite sequence $\left\{ \mb Z_{i}\right\} $ of independent,
random, self-adjoint matrices. For all $t\in\mbb R$, 
\begin{align*}
\P\left(\lambda_{\max}\left(\sum_{i=1}^{n}\mb Z_{i}\right)\ge t\right) & \le\inf_{\theta>0}e^{-\theta t}\tr\exp\left(\sum_{i=1}^{n}\log\E e^{\theta\mb Z_{i}}\right).
\end{align*}
\end{thm}

 In particular, we combine Lemma \ref{lem:KS} and Theorem \ref{thm:master-tail}
to obtain a sharper analog to the tail bounds for Rademacher series
derived in \citep{tropp_user-friendly_2011}, for general centered
Bernoulli random variables. As a consequence of the use of Kearns-Saul inequality (i.e., Lemma \ref{lem:KS}), the improvement over similar bounds obtained via matrix Hoeffding or matrix Bernstein inequalities can be particularly significant if the Bernoulli random variables have means close to $0$ or $1$.

\begin{prop}
\label{prp:Bernoulli-series-tail} For $i=1,2,\dots,n$ let $\delta_{i}\sim\mr{Bernoulli}(p_{i})$
be independent random variables. Furthermore, let $\mb X_{i}$ be
deterministic $N\times N$ self-adjoint matrices. Then with $K_{i}=K(p_{i})$
defined as in Lemma \ref{lem:KS}, we have 
\begin{align*}
\P\left(\left\lVert \sum_{i=1}^{n}\left(\delta_{i}-p_{i}\right)\mb X_{i}\right\rVert \ge t\right) & \le2Ne^{-\frac{t^{2}}{4\sigma^{2}}},
\end{align*}
 where 
\begin{align*}
\sigma^{2} & =\left\lVert \sum_{i=1}^{n}K_{i}^{2}\mb X_{i}^{2}\right\rVert .
\end{align*}
\end{prop}
\begin{IEEEproof}
Let $\mb Z_{i}=\left(\delta_{i}-p_{i}\right)\mb X_{i}$ for $i=1,2,\dotsc,n$.
For any real number $\theta$ we have 
\begin{align*}
\E e^{\theta\mb Z_{i}} & =p_{i}e^{\theta\left(1-p_{i}\right)\mb X_{i}}+(1-p_{i})e^{-\theta p_{i}\mb X_{i}}.
\end{align*}
 Since $\theta\mb X_{i}$ is self-adjoint, it can be diagonalized.
Therefore, by applying Lemma \ref{lem:KS} to the eigenvalues of $\theta\mb X_{i}$
the above equation implies that 
\begin{align*}
\E e^{\theta\mb Z_{i}} & \preccurlyeq e^{\theta^{2}K_{i}^{2}\mb X_{i}^{2}},
\end{align*}
 where the inequality is with respect to the positive semidefinite
cone. Therefore, we have 
\begin{align*}
\tr\exp\left(\sum_{i=1}^{n}\log\E e^{\theta\mb Z_{i}}\right) & \le\tr\exp\left(\theta^{2}\sum_{i=1}^{n}K_{i}^{2}\mb X_{i}^{2}\right)\\
 & \le\! N\,\exp\left(\theta^{2}\left\lVert \sum_{i=1}^{n}K_{i}^{2}\mb X_{i}^{2}\right\rVert \right)=\!N e^{\theta^{2}\sigma^{2}}.
\end{align*}
 Then it follows from Theorem \ref{thm:master-tail} that 
\begin{align*}
\P\left(\lambda_{\max}\left(\sum_{i=1}^{n}\mb Z_{i}\right)\ge t\right) & \le\inf_{\theta>0}N\,e^{-\theta t+\theta^{2}\sigma^{2}}=Ne^{-\frac{t^{2}}{4\sigma^{2}}}.
\end{align*}
 Replacing $\mb X_{i}$ by $-\mb X_{i}$ and repeating the above argument
we can similarly show that 
\begin{align*}
\P\left(\lambda_{\min}\left(\sum_{i=1}^{n}\mb Z_{i}\right)\le-t\right) & \le Ne^{-\frac{t^{2}}{4\sigma^{2}}}.
\end{align*}
The union bound then guarantees that 
\begin{align*}
\P\left(\left\lVert \sum_{i=1}^{n}\mb Z_{i}\right\rVert \ge t\right) & \le2Ne^{-\frac{t^{2}}{4\sigma^{2}}},
\end{align*}
as desired.
\end{IEEEproof}

\section*{Acknowledgements}
\addcontentsline{toc}{section}{Acknowledgment}
S. Bahmani and J. Romberg were supported in part by ONR grant N00014-11-1-0459, NSF grants CCF-1415498 and CCF-1422540, and the Packard Foundation. P. Tetali was supported in part by NSF grants DMS-1407657 and CCF-1415498.

\bibliographystyle{abbrvnat}
\bibliography{references}

\begin{thebibliography}{16}
\providecommand{\natexlab}[1]{#1}
\providecommand{\url}[1]{\texttt{#1}}
\expandafter\ifx\csname urlstyle\endcsname\relax
  \providecommand{\doi}[1]{doi: #1}\else
  \providecommand{\doi}{doi: \begingroup \urlstyle{rm}\Url}\fi

\bibitem[Ben-Shimon and Krivelevich(2009)]{ben-shimon_vertex_2009}
S.~Ben-Shimon and M.~Krivelevich.
\newblock Vertex percolation on expander graphs.
\newblock \emph{European Journal of Combinatorics}, 30\penalty0 (2):\penalty0
  339--350, 2009.

\bibitem[Chung and Horn(2007)]{chung_spectral_2007}
F.~Chung and P.~Horn.
\newblock The spectral gap of a random subgraph of a graph.
\newblock \emph{Internet Math.}, 4\penalty0 (2-3):\penalty0 225--244, 2007.

\bibitem[Chung et~al.(2009)Chung, Horn, and Lu]{chung_giant_2009}
F.~Chung, P.~Horn, and L.~Lu.
\newblock The giant component in a random subgraph of a given graph.
\newblock In \emph{International Workshop on Algorithms and Models for the
  Web-Graph}, pages 38--49. Springer, 2009.

\bibitem[Colbourn(1987)]{colbourn_combinatorics_1987}
C.~J. Colbourn.
\newblock \emph{The Combinatorics of Network Reliability}, volume~4 of
  \emph{International Series of Monographs on Computer Science}.
\newblock Oxford University Press, Inc., New York, 1987.

\bibitem[Greenhill et~al.(2008)Greenhill, Holt, and
  Wormald]{greenhill_expansion_2008}
C.~Greenhill, F.~B. Holt, and N.~Wormald.
\newblock Expansion properties of a random regular graph after random vertex
  deletions.
\newblock \emph{European Journal of Combinatorics}, 29\penalty0 (5):\penalty0
  1139--1150, 2008.

\bibitem[Grimmett(1999)]{grimmett_percolation_1999}
G.~Grimmett.
\newblock \emph{Percolation}, volume 321 of \emph{Grundlehren der
  mathematischen {Wissenschaften}}.
\newblock Springer Berlin Heidelberg, Berlin, Heidelberg, 1999.

\bibitem[Harris and Srinivasan(2014)]{harris_improved_2014}
D.~G. Harris and A.~Srinivasan.
\newblock Improved bounds and algorithms for graph cuts and network
  reliability.
\newblock In \emph{Proceedings of the Twenty-Fifth Annual {ACM-SIAM} Symposium
  on Discrete Algorithms}, {SODA} '14, pages 259--278, Philadelphia, PA, USA,
  2014. {SIAM}.

\bibitem[Karger(2001)]{karger_randomized_2001}
D.~R. Karger.
\newblock A randomized fully polynomial time approximation scheme for the
  all-terminal network reliability problem.
\newblock \emph{{SIAM} Review}, 43\penalty0 (3):\penalty0 499--522, 2001.

\bibitem[Karger(2016)]{karger_fast_2016}
D.~R. Karger.
\newblock A fast and simple unbiased estimator for network (un)reliability.
\newblock In \emph{Foundations of Computer Science ({FOCS}), 2016 {IEEE} 57th
  Annual Symposium on}, 2016.

\bibitem[Kearns and Saul(1998)]{kearns_large_98}
M.~Kearns and L.~Saul.
\newblock Large deviation methods for approximate probabilistic inference.
\newblock In \emph{Proceedings of the Fourteenth Conference on Uncertainty in
  Artificial Intelligence}, UAI'98, pages 311--319, San Francisco, CA, USA,
  1998. Morgan Kaufmann Publishers Inc.
\newblock ISBN 1-55860-555-X.

\bibitem[Krivelevich(2016)]{krivelevich_phase_2016}
M.~Krivelevich.
\newblock The phase transition in site percolation on pseudo-random graphs.
\newblock \emph{The Electronic Journal of Combinatorics}, 23\penalty0
  (1):\penalty0 1--12, 2016.

\bibitem[Levin et~al.(2009)Levin, Peres, and Wilmer]{levin_markov_2009}
D.~A. Levin, Y.~Peres, and E.~L. Wilmer.
\newblock \emph{Markov chains and mixing times}.
\newblock American Mathematical Society, 2009.

\bibitem[Mohar(1989)]{mohar_isoperimetric_1989}
B.~Mohar.
\newblock Isoperimetric numbers of graphs.
\newblock \emph{Journal of Combinatorial Theory, Series B}, 47\penalty0
  (3):\penalty0 274--291, 1989.

\bibitem[Oliveira(2009)]{oliveira_concentration_2009}
R.~I. Oliveira.
\newblock Concentration of the adjacency matrix and of the {Laplacian} in
  random graphs with independent edges.
\newblock \texttt{arXiv:0911.0600 [math.CO]}, Nov. 2009.

\bibitem[Sun et~al.(2006)Sun, Boyd, Xiao, and Diaconis]{sun_fastest_2006}
J.~Sun, S.~Boyd, L.~Xiao, and P.~Diaconis.
\newblock The fastest mixing {Markov} process on a graph and a connection to a
  maximum variance unfolding problem.
\newblock \emph{{SIAM} review}, 48\penalty0 (4):\penalty0 681--699, 2006.

\bibitem[Tropp(2011)]{tropp_user-friendly_2011}
J.~A. Tropp.
\newblock User-friendly tail bounds for sums of random matrices.
\newblock \emph{Foundations of Computational Mathematics}, Aug. 2011.

\end{thebibliography}

\end{document}